\documentclass[11pt]{article}
\usepackage{latexsym}
\usepackage{amsmath,amsthm,amssymb}

\newtheorem{theorem}{Theorem}[section]
\newtheorem{proposition}[theorem]{Proposition}
\newtheorem{lemma}[theorem]{Lemma}
\newtheorem{corollary}[theorem]{Corollary}

\title{Index theory with bounded geometry, the uniformly finite $\hat{A}$ class, 
and infinite connected sums}
\author{Kevin Whyte}
\date{}
\begin{document}

\maketitle
\section*{Introduction}

We work primarily in the category of manifolds of bounded geometry.  The
objects are manifolds with bounds on the curvature tensor, its
derivatives, and on the injectivity radius.  The morphisms are
diffeomorphisms of bounded distortion.  We think of these manifolds as
having  a chosen bounded distortion class of metrics.   Unless otherwise
stated, all manifolds in the paper are assumed of  this type.  

These definitions are designed to reflect the restrictions 
imposed on a non-compact manifold which is controlled in some way by  
a compact manifold.  The most common example of this is a covering of 
a compact manifold.  Any metric on the base gives a metric of bounded 
geometry on the cover, and any two such metrics are bounded 
distortion equivalent.  Similarly, leaves of foliations of compact 
manifolds have a canonical bounded distortion class of metric of 
bounded geometry.

We try to understand index theory for these manifolds, and in 
particular, questions of positive scalar curvature.  Generally, the 
appropriate notion here is uniformly positive scalar curvature, 
meaning the scalar curvature is bounded away from zero from below.   
When we say that $M$ admits a metric of positive scalar curvature, we 
mean that within the chosen bounded distortion class of metrics there 
is a metric of uniformly positive scalar curvature. 

To understand positive scalar curvature we need an appropriate 
generalization of the $\hat{A}$ class.  One interesting feature which 
emerges is that this class lives in a non-Hausdorff homology group, 
and thus standard $C^{*}$ algebra methods do not apply.  It turns out 
that to understand this class requires rather delicate spectral 
estimates for Dirac operators on the boundaries of certain compact 
submanifolds.  The resulting theorems have unexpected applications to 
compact manifolds.

The motivation for this work comes from some interesting infinite 
connected sum examples studied in \cite{BW1} and \cite{Roe}.  As these 
examples provide good motivation for the definition of the $\hat{A}$ 
class we use, and are interesting on their own, we discuss them in 
some detail in the first section.  The gap between the obstruction of 
\cite{Roe} and the construction in \cite{BW2} is a direct consequence 
of the non-Hausdorff problem mentioned above.  Our theorems solve 
this problem, and provide a complete picture of these examples.  

We would like to thank Shmuel Weinberger for bringing this problem to 
our attention, and for many helpful discussions.  We would also like to
thank the referee, without whom this paper would be substantially less
readable.

\section{Infinite Connected Sums}

  Let $M$ be a manifold of bounded geometry 
and $S$ a discrete subset of $M$.  Given $N$, a compact manifold, we 
form a new manifold, $M\#_{S}N$, by connected summing a copy of $N$ at 
each point of $S$.  For this manifold to have a well defined bounded 
distortion class of metrics of bounded geometry one needs that points 
of $S$ are uniformly separated, meaning that for some $\varepsilon > 
0$ any two distinct points of $S$ are at least $\varepsilon$ apart.  
Such a set is called {\em uniformly discrete}.  One important example 
of this construction is for $M$ the universal cover of some $X$ and 
$S$ an orbit of $\pi_{1}(X)$.  Then $M\#_{S}N$ is just a cover of $X\#N$.

  Now suppose $M$ has positive scalar curvature.  We wish to understand 
when $M\#_{S}N$ admits a metric of positive scalar curvature.  We 
restrict our attention to the case of $M$ and $N$ spin, and $N$ simply 
connected.  If we are in the covering space situation, there is a positive 
scalar curvature metric invariant under the covering group if and only if 
$\hat{A}(N)=0$.  We are looking to generalize the $\hat{A}$ class as 
an obstruction to positive scalar curvature, so we will typically 
assume $\hat{A}(N) \neq 0$.

  A natural first question is whether it is ever possible, keeping 
bounded geometry, to have a metric of positive scalar curvature in 
the presence of any such ``obstructing'' $N$.  Examples showing this is
possible are constructed in \cite{BW2}.  We sketch their construction 
below as it gives insight about the kind of obstructions that arise 
in our main theorem, Theorem \ref{poscal}.

\begin{lemma} (\cite{BW2})  Let $M$ be non-compact, spin, and of 
uniformly positive scalar curvature, and $N$ be compact, spin, and 
simply connected.  The manifold $M\#N$  admits a metric of uniformly 
positive scalar curvature.
\begin{proof}

   Choose a geodesic ray $p$ in $M$.  Let $\hat{M}$ be the infinite 
 connected sum of $M$ with copies of $N$ at $p(2i)$ and 
 copies of $\bar{N}$ ($N$ with opposite orientation) at $p(2i+1)$, for every 
 $i$.  We view $\hat{M}$ in two ways.  First, we think of the copies 
 of $N$ as coming in pairs, at $p(2i)$ and $p(2i+1)$.  Each of these 
 pairs is a copy of $N \# \bar{N}$.  This manifold is null cobordant, 
 indeed it is the boundary of $N \times [0,1]$.  Further, as $N$ is 
 simply-connected, the null cobordism can be realized by a sequence 
 of surgeries of index greater than $1$.  
 
   The collection of all these surgeries is uniformly locally finite, 
and thus can be carried out simultaneously without leaving the bounded 
geometry category.  Using surgery in positive scalar curvature (see \cite{GL}
 and \cite{SY}), we can also carry out the surgeries keeping positive scalar 
curvature. Note that to keep the metric in bounded distortion class of the 
connected sum, one changes the metric only in a neighborhood of 
the sphere one is surgering.  The resulting distortion in that neighborhood
depends on the lower bound for the scalar curvature in $M$, and thus
the assumption of uniformly positive scalar curvature cannot be weakened
to merely positive scalar curvature.

   Now view $\hat{M}$ similarly, but with the copies of $N$ paired 
as $p(2i-1)$ and $p(2i)$.  One can similarly carry out surgeries here, 
leaving just a single copy of $N$ at $p(0)$.  Thus, since $\hat{M}$ 
has positive scalar curvature so does $M\#N$. 
  
\end{proof}
\end{lemma}

  Notice that we did not really need that $p$ was a geodesic, just 
that the time it spends in any ball is uniformly bounded in the 
radius.  Thus to carry out the above construction for $M\#_{S}N$, one 
needs such tails from each point of $S$ which spread out in such a 
way that only a uniformly bounded number pass through any ball.  The 
existence of such ``tails'' is a homological question.

  We say that $c$, an $i$-chain in $M$, is uniformly finite if there 
is a bound on the diameter of simplices in the support of $c$ and for 
every $r$ there is an upper bound $C_{r}$ on the sum of the absolute 
values of the coefficients of the simplices which intersect any $r$-ball.  
We denote these chains by $C_{i}^{uf}(M)$, and the corresponding homologies 
by $H_{i}^{uf}(M)$.  

  Any $S$ as above gives a natural element of $C_{0}^{uf}(M)$, and a 
collection of tails gives a null homology of that class.  We denote 
that class corresponding to $S$ in $H_{0}^{uf}(M)$ by $[S]$.  We have 
sketched the proof of:

\begin{theorem}\label{constr}(\cite{BW2}) Let $M$, $N$, and $S$ be as above.  If 
$[S]=0$ then $M\#_{S}N$ admits a metric of positive scalar curvature.
\end{theorem}

  To get a feel for the nature of the obstruction measured by 
$H_{0}^{uf}$, consider the lattice ${\mathbb Z}^{2}$ in ${\mathbb 
R}^{2}$.  For any $r$, the number of lattice points in a ball of 
radius $r$ is about $r^{2}$, while the perimeter is about $r$.  Thus, 
as $r$ increases, the number of tails crossing the boundary in some 
bounded region must increase unboundedly.  This shows that $[{\mathbb 
Z}^{2}]$ is non-zero.

  This is essentially the only obstruction to finding tails.
  
\begin{theorem} \label{growth}(\cite{BW1}, \cite{Why})  If $c \in C_{0}^{uf}(M)$ 
then $[c]=0$ if and only if there are $r$ and $C$ such that for any $R \subset M$ 
one has:

$$|\Sigma_{\sigma \in R} c_{\sigma}| \leq C vol(\partial_{r} R)$$  
\end{theorem}

  Recall that a {\em regular sequence} is a sequence $R_{i}$ of subsets 
of $M$ for which for any $r$, 
$$\lim_{i \to \infty} \frac{vol(\partial_{r}R_{i})}{vol(R_{i})}=0$$

 $M$ is called {\em amenable} if and only if there is a regular sequence in $M$.

\begin{theorem} \label{nonamen}(\cite{BW2}) The following are equivalent:
\begin{enumerate}
\item $M$ is non-amenable.
\item $H_{0}^{uf}(M)=0$.
\item For all $S$ uniformly discrete, $[S]=0$.
\end{enumerate}
\end{theorem}

 This means, in particular, that if, in the context of Theorem \ref{constr},
$M$ is non-amenable then for any $N$ and $S$, $M\#_{S}N$ admits a metric of positive scalar curvature.  If $M$ is 
a universal cover of a compact manifold then amenability is equivalent
to amenability of the fundamental group. Thus, for example, let $X$ be a
surface of higher  genus cross $S^{4}$ connected sum any $N$.  If
$\hat{A}(N)\neq 0$ then 
$X$ does not admit a metric of positive scalar curvature, but the 
universal cover does have a metric of uniformly positive scalar curvature
in the bounded distortion class of the periodic metrics. 

  Naturally one wants to know whether $[S]$ in $H_{0}^{uf}$ really is 
an obstruction to a metric of positive scalar curvature.  That it is 
follows from the theory of the next section.

\section{The uniformly finite $\hat{A}$ genus}

  How does one obstruct metrics of positive scalar curvature?  In the 
compact case one has the $\hat{A}$ genus.  We want to generalize this 
to the bounded geometry setting.

  According to Chern-Weil theory (see, for example, \cite{LM}), the 
$\hat{A}$ class can be defined as the cohomology class of a universal 
polynomial, $P$, in the curvature tensor.  On a manifold of bounded
geometry,  this form is bounded, and thus represents an element of
$l^\infty$-cohomology.  This cohomology class depends only on the bounded
distortion diffeomorphism class of the metric.  The proof on this is
similar to the standard proof of invariance of characteristic classes:

Let $\{g_t\}$ be a one parameter family of metrics on $M$ for which the
induced metric $g_t + {dt}^2$ on $M \times [0,1]$ has bounded
curvature.   Let $\omega$ be the characteristic form, given by $P$, on $M
\times [0,1]$.  Write $\omega=\alpha_t + \beta_t \wedge {dt}$, where
$\alpha_t$ and $\beta_t$ are forms on $M$.  By the assumption of bounded
curvature, these forms are bounded.   Since $\omega$ is closed,
$\frac{d}{dt} \alpha_t = d \beta_t$.  Thus $\alpha_1 - \alpha _0 =
d (\int_0^1 \beta_t)$, and therefore $\alpha_0$ and $\alpha_1$ define
the same $l^\infty$-cohomology class on $M$.  As these are the
characteristic forms on $M$ for $g_0$ and $g_1$, we see that the
$l^\infty$-cohomology class of the characteristic form is the same for
metrics connected by such one parameter families.

For any two metrics of bounded curvature on $M$, $g_0$ and $g_1$,
which are bounded distortion equivalent, lemma 2.6 of \cite{CG2} shows
that the one parameter family of metric $g_t = tg_1 + (1-t)g_0$ has
bounded curvature, provided that the difference of Levi-Civita
connections is a bounded operator (note that by the assumption of bounded
distortion equivalence, bounded means the same thing for all the
metrics).  It is easy to see that this is the case if the identity map
$(M,g_0)$ to $(M,g_1)$ is not only bounded distortion, but also has
bounded $2$-jet.  

Let $g_0$ and $g_1$ be any two bounded distortion equivalent metrics of
bounded geometry on $M$.  By Theorem 2.5 of \cite{CG2}, we may assume
that not only do these metrics have bounded curvature, but that their
curvature tensors have bounded covariant derivatives to arbitrary order. 
Thus there is some $r>0$ so that the exponential maps on the balls of
radius $r$ are bounded distortion diffeomorphisms with bounded
derivatives of arbitrary order.   The standard proof that $C^1$
diffeomorphic smooth manifolds are $C^\infty$ diffeomorphic by
convolution with a smoothing kernel (see, for example, \cite{hirsch},
section 2.2), shows that there is a map $f:(M,g_0) \to (M,g_1)$ at
finite distance from the identity, which is a bounded distortion
diffeomorphism with bounded 2-jet.  Since $f$ induces the identity on
top dimensional $l^\infty$ cohomology, this shows that there is a well
defined $l^\infty \hat{A}$ genus for every bounded distortion
diffeomorphism class of bounded geometry metrics on $M$. 

  This fits with the previous discussion as $l^{\infty}$ cohomology 
is naturally Poincare dual to uniformly finite homology (\cite{AW}).  As
the proof is  simple in the case we use, we include it here:

\begin{lemma} Let $M^m$ be a complete, connected, Riemannian manifold of
bounded geometry.  There is a canonical isomorphism between $H_0^{uf}(M)$
and $H^m_\infty(M)$.
\end{lemma}

\begin{proof}

  Let $\varepsilon > 0$ be much smaller than the convexity radius of
$M$.  Let $S$ be a maximal subset of $M$ such that any two points of $S$
are at distance at least $\varepsilon$.  The balls of radius
$\varepsilon$ centered at points of $S$ cover $M$, and the concentric
balls of radii $\frac{\varepsilon}{2}$ are disjoint.  As $S$ with its
induced metric is quasi-isometric to $M$, $H_0^{uf}(S)=H_0^{uf}(M)$.

 Choose a partition of unity $\{f_s\}_{s \in S}$ so that $f_s$ is
supported in $B_s(\varepsilon)$.   Given the bounds on the geometry of
$M$ these can be chosen with uniformly bounded derivatives.  Let $\phi_s$
be the bump form $\frac{f_s dvol}{\int f_s}$.  

Given $c \in C_0^{uf}(S)$, let $w_c = \Sigma_s c_s \phi_s$.  If
$c=\partial b$ for $b \in C_1^{uf}(S)$ then $b$ is a uniformly locally
finite sum of pairs $(s,s')$ with $d(s,s')$ uniformly bounded.  The
difference $\phi_s - \phi_{s'}$ is therefore $d$ of a  bounded $n-1$ form
of uniformly bounded support. Thus $[w_c]=0$  in $l^\infty$ cohomology,
so $c \mapsto w_c$ induces a well defined map $H_0^{uf}(S) \to
H^n_{\infty}(M)$.

Similarly, given $w$ an $l^\infty$ $n$-form on $M$, define $c_w=\Sigma_s
(\int f_s w) s$.  If $w=d\eta$ for an $l^\infty$ form $\eta$, then by
Theorem \ref{growth} and Stokes' Theorem, $c_w=0$ in $H_0^{uf}(S)$. 
Thus we have a map $H^n_{\infty}(M) \to H_0^{uf}(S)$.  We now show these
maps are inverses.

Given a uniformly finite chain $c=\Sigma c_s s$ let  $c'$ be the image of
$c$ under the composition of these maps.  We have $c'=\Sigma_s c_s
\Sigma_t (\int f_t \phi_s) t$.  Let $d_s = \Sigma_t (\int f_t \phi_s)
t$, then $c-c' = \Sigma_s c_s ( s - d_s)$.  The chain $s - d_s$ is
of uniformly bounded support and sum to zero, hence it is a boundary
of a uniformly bounded $1$-chain, $b_s$ with support in a uniformly bounded
neighborhood of $s$.   Thus the difference $c-c' = \partial b = \partial
\Sigma_s c_s b_s$ with $b \in C_1^{uf}$.  

Likewise, given an $l^\infty$ form $w$, let $w'$ be the image under
the composition of the two maps.  We have $w=\Sigma_s f_s w$, and
$w'=\Sigma_s (\int f_s w) \phi_s$.  Thus $w-w'= \Sigma_s \eta_s$, where  
the forms $\eta_s =  f_s w - (\int f_s w) \phi_s$ are of uniformly bounded
support, uniformly bounded pointwise norm, and have integral zero.  This
implies, by the bounded geometry of $M$, that $\eta_s = d \sigma_s$ for
$\sigma_s$ also of uniformly bounded norm and support.  Thus $w-w'=
d(\Sigma_s \sigma_s)$, and therefore $w=w'$ in $H^n_\infty(M)$.

Thus the maps are inverses and give the desired isomorphism. 
\end{proof}

We define the uniformly finite $\hat{A}$ genus, $\hat{A}^{uf}$, as the
dual in $H_{0}^{uf}$  of the $\hat{A}$ form.  For the infinite connected
sums we have been discussing, one has
$\hat{A}^{uf}(M\#_{S}N)=\hat{A}^{uf}(M) +
\hat{A}(N)[S]$,  where $\hat{A}(N)$ is the (integer) $\hat{A}$-genus of
$N$.

\begin{theorem}\label{poscal} If $M$ has non-negative scalar curvature 
then $\hat{A}^{uf}=0$.
\end{theorem}

  The proof of this theorem is, in outline, much like the 
corresponding theorem in the compact case: relate the $\hat{A}$ class 
to the index of a Dirac operator via the index theorem and then prove 
the vanishing of this index by a Bochner type argument.  Both steps 
are substantially more difficult here, and some essentially new 
ingredients are needed.

  We first need to reinterpret Theorem \ref{growth} cohomologically.
Cheeger  and Gromov (\cite{CG1}) prove a chopping theorem for manifolds
of bounded  geometry which says that for any $n$ there are constants 
$C_0,C_1,...,C_n$ and $r$ so that for any $S\subset M$ there is a codimension 
$0$ manifold with boundary $(X,\partial X)$ such that:

\begin{enumerate}
\item $S \subset X \subset N_r(S)$
\item For $i=0,1,...,n$ $\nabla^i {II_{\partial X}} \leq C_i$ where 
$II_{\partial X}$ is the second fundamental form, and $\nabla^i$ is 
the $i^{th}$ covariant derivative.
\item $\frac{Vol(\partial X)}{Vol(\partial_r(S))}$ is bounded above and 
below independent of $S$.
\end{enumerate}

\begin{lemma} \label{cohvan}For any $n$ there are constants 
$C_0,C_1,...,C_n$ and $r$ so that if $\omega \in \Omega^m(M)$ is
bounded then $\omega$ is $d$ of a bounded form if and only if for some $C$

$$|\int_X \omega| \leq C vol(\partial X)$$

for all $(X,\partial X)$ compact, codimension 0 submanifold with 
$\nabla^i {II_{\partial X}} \leq C_i$ for $i=0,1,...,n$.
\begin{proof}

The condition is necessary by Stokes theorem.  To see that it is 
sufficient we use Theorem \ref{growth}.  That lemma, formulated
cohomologically, says that we need to show: 

$$|\int_S \omega| \leq C vol(\partial_r S)$$

With an arbitrary $S \subset M$ in place of $X$.  

For any such $S$, approximate it by a manifold with boundary, $X$,
via the chopping theorem above.  

By (1) $|\int_S \omega|$ and $|\int_X \omega|$ differ by at most 
$||\omega|| Vol(\partial_r S)$ which, by part (3) of the chopping 
theorem, is bounded above by $K vol(\partial X)$ for some $K$
which depends on $\omega$ and $M$, but not on $S$.

Thus the bound on $|\int_S \omega|$ for arbitrary $S$ follows, with 
perhaps a larger $C$, from the bound for submanifolds with bounded 
second fundamental forms. 

\end{proof} 
\end{lemma}

  In view of Lemma \ref{cohvan}, Theorem \ref{poscal} will follow 
from:

\begin{theorem}\label{bndry} Let $(X,\partial X)$ be a compact spin 
manifold of non-negative scalar curvature.  There is a $C$ depending only 
on the curvature and second fundamental form so that
$$|\int_{X}\hat{A}|\leq C vol(\partial X)$$
\end{theorem}

  We prove Theorem \ref{bndry} by a detailed analysis of the Dirac 
operator on a manifold with boundary.

\section{Index Theory}

Let $D$ be the canonical Dirac operator on spinors (much of this 
section works for an arbitrary geometric operators, but we will not 
need this generality).  
  
\begin{theorem}(\cite{APS}) For any geometric operator $D$, 
there is a characteristic form $\omega$, a polynomial $P$, and $n \in 
{\mathbb N}$, so that for any $(X,\partial X)$ compact we have:

$$index(D) = \int_X w + \int_{\partial X} P(II_{\partial X},\nabla 
{II}_{\partial X},...,\nabla^n {II}_{\partial X}) + \eta(\partial 
X)$$ \label{aps}

Where $\eta(N)=\lim_{s\rightarrow 0^+} \Sigma_{\lambda \in 
Spec(D_\partial)} \lambda^{-s}$ 

\end{theorem}

  We will use this to prove vanishing of $\omega$ in 
$l^{\infty}$-cohomology via Lemma \ref{cohvan}. 

  Given the bounds on $\nabla^i {II}_{\partial X}$, the middle term in 
Theorem \ref{aps} is bounded by a constant multiple of $Vol(\partial X)$.

Likewise, $\eta$ is bounded linearly in $Vol(\partial X)$.  This 
is shown for the signature operator in \cite{CG2} and for a wide 
range of geometric operators including the Dirac operator on spinors 
in \cite{Ra}. 

In view of this, and Lemma \ref{cohvan}, we can interpret \ref{aps} as saying that the 
characteristic form in $H_0^{uf}$ is a ``uniformly finite index'' of our operator.    

  For the Dirac operator on spinors, the form $\omega$ is the 
$\hat{A}$ form.  The Dirac operator is related to positive scalar 
curvature by:

\begin{theorem}\label{lic}(Lichnerowicz formula, \cite{LM})

  $$D^2 = \nabla^{*} \nabla + \frac{\kappa}{4}$$

Where $\nabla$ is the canonical spinor connection, and $\kappa$ the 
scalar curvature.

\end{theorem}
 
So, if $s$ is a harmonic spinor (meaning $Ds=0$) on $X$, we have:

$$0 = <\nabla^{*} \nabla s, s> + \frac{\kappa}{4} ||s||^2 $$ 

If $X$ were closed, we could integrate over $X$, and the first term 
on the RHS would be $||\nabla s||^2$.

Then, if $\kappa>0$ both terms on the RHS would be $\geq 0$, and 
therefore 0.  This would mean $s=0$. i.e. that there are no harmonic 
spinors, so that the index would be zero.  Then by the Atiyah-Singer 
index theorem the $\hat{A}$ genus would be 0. 

It is this argument we try to extend to our setting.

When one integrates the Lichnerowicz formula over an $X$ with 
boundary there is an extra term which comes from the boundary term of 
integration by parts (see \cite{LM}).

$$\int_X <\nabla^{*} \nabla s,s> = ||\nabla s||^2  - \int_{\partial 
X} <\nabla_\nu s,s>$$
($\nu$ is the unit normal vector to $\partial X$): 
 
The second term on the RHS is introduced when we integrate by parts. 
     
So, when we have boundary, the Lichnerowicz formula becomes:

\begin{lemma} If $s$ is a harmonic spinor on $X$, we have:

 $$0 = ||\nabla s||^2 + \int_X \frac{\kappa}{4} ||s||^2 - 
\int_{\partial X} <\nabla_\nu s,s>$$

\end{lemma}

If we assume the scalar curvature is $\geq 0$ then the first two 
terms are non-negative.  This can only happen if:

$$\int_{\partial X} <\nabla_\nu s,s> \geq 0$$

We can expand $D$ in normal coordinates around the boundary:

\begin{lemma}(\cite{Gi}) Along $\partial X$ we have: 

$$D = G ( D_{\partial X} - \nabla_\nu - \frac{1}{2} tr(II))$$

Where $G$ is the bundle automorphism induced by clifford 
multiplication by the normal vector, and $D_{\partial X}$ is the 
Dirac operator intrinsic to $\partial_X$. 

\end{lemma}
 
$Ds=0$ gives

$$<\nabla_\nu s,s> = <D_{\partial X}s,s> - \frac{1}{2} tr(II) 
||s||^2$$

Thus we have proven:

\begin{proposition}\label{void} If $X$ has non-negative scalar curvature, 
and $s$ is a harmonic spinor, then:

$$\int_{\partial X} <D_{\partial X}s,s> \geq \int_{\partial X} 
\frac{1}{2} tr(II) ||s||^2$$

\end{proposition}

The boundary conditions for harmonic spinors in the 
Atiyah-Patodi-Singer index theorem are that when $s|_{\partial X}$ is 
expanded in eigenfunctions of $D_{\partial X}$, only negative 
eigenvalues are used.  

Writing $s$ as $\Sigma_{\lambda} a_\lambda s_\lambda$, where the 
$s_\lambda$ are the eigenvectors of $D_{\partial X}$, \ref{void} 
becomes:

$$ \Sigma_\lambda {\lambda ||a_\lambda||^2} \geq \frac{1}{2} 
\int_{\partial X} tr(II) ||s||^2 $$

  Since all the $\lambda$ must be negative, we must have some 
$\lambda \geq \lambda_0={inf}(\frac{1}{2} tr(II))$.  By projecting 
onto these eigenspaces between $\lambda_0$ and $0$, we get:

\begin{theorem} If $(X,\partial X)$ is spin and has non-negative scalar 
curvature then there is $\Lambda$, depending only on the second 
fundamental form, such that $dim(H) \leq N_{D^2_{\partial 
X}}(\Lambda)$, where $H$ is the space of harmonic spinors on $X$ with 
A-P-S boundary values, and $N_{D^2_{\partial X}}(\Lambda)$ is the 
dimension of the space of eigenfunctions of $D^2_{\partial X}$ below 
$\Lambda$.
\end{theorem}

\begin{theorem}\label{spec} If $N^n$ is a compact spin manifold, then for 
each $\lambda$ there is a $C_\lambda$ depending only on the curvature 
and injectivity radius of $N$, for which $N_{D^2}(\lambda) \leq 
C_\lambda vol(N)$

\begin{proof}

By (\cite{Bo},Prop 4.20(ii)) there are $A$ and $B$ so that 
$\lambda_n(D^2) \geq \lambda_{An-B}(\Delta)$, where $\Delta$ is the 
laplacian on functions.  Thus the bound we need follows immediately 
from the same statement (with different $C_\lambda$) for the 
laplacian on functions.

\begin{theorem}(\cite{Gr},Appendix $C_+$) There is a constant $K$
depending only on  curvature bounds and dimension for which for any
compact manifold $N$  we have the bound:

$$ \lambda_{V(\varepsilon)} \geq K\varepsilon^{-2} $$

For any $\varepsilon \leq inj\ rad$.  Here $V(\varepsilon)$ is the 
minimal number of $\varepsilon$ balls which cover $N$.

\end{theorem}

Since we have bounds on the curvature, there is an $L$ such that:

$$Vol(B(\varepsilon)) \geq L\varepsilon^n $$ 

for $B(\varepsilon)$ any such $\varepsilon$ ball in $N$.   

Choose a maximal family of disjoint $\frac{\varepsilon}{2}$ ball in 
$N$.  By maximality the concentric balls of radius $\varepsilon$ 
cover.  But by disjointness there are at most 

$$\frac{2^n Vol(N) \varepsilon^{-n}}{L}$$

balls.  Therefore we have:

\begin{proposition} There is a constant $C$ depending only on the curvature 
such that, for any $\varepsilon \leq inj\ rad (N)$ 

$$\lambda_{C vol(N) \varepsilon^{-n}} \geq \epsilon^{-2} $$
\end{proposition}
 
Turning this around to an upper bound on the spectral counting 
function gives: 

\begin{proposition} There is a constant $C$ depending only on the curvature, 
and a $\lambda_0$ depending only on curvature and injectivity radius, 
such that for any $\lambda \geq \lambda_0$ 

$$N(\lambda) \leq C Vol(N) \lambda^\frac{n}{2}$$
\end{proposition}

This gives the $C_\lambda$ we needed for \ref{spec}
 
$C_\lambda = C \lambda^\frac{n}{2}$ works for $\lambda \geq 
\lambda_0$.  

and  

$C_\lambda=C \lambda_0^\frac{n}{2}$ 

works for $\lambda \leq \lambda_0$, as 

$N(\lambda) \leq N(\lambda_0)$

\end{proof}
\end{theorem}

This completes the proof of Theorem \ref{bndry}, and thereby Theorem 
\ref{poscal}.

\section{Applications}

Our first corollary, combined with the work of \cite{BW2}, gives a 
complete characterization of infinite connected sums of positive 
scalar curvature.

\begin{corollary}\label{sums} Let $M$, $N$, and $S$ be as before.  If 
$\hat{A}(N)\neq 0$ then $M\#_{S}N$ admits a metric of uniformly positive 
scalar curvature if and only if $[S]=0$ in $H_{0}^{uf}(M)$.
\end{corollary}

  In fact, \ref{poscal} shows that if $[S]\neq 0$ then $M\#_{S}N$ does 
not even admit a metric of non-negative scalar curvature.  Thus there 
is an interesting alternative: an infinite connected sum where $M$ has 
uniformly positive scalar curvature either has a metric of uniformly 
positive scalar curvature or does not even have a metric of 
non-negative scalar curvature.  It would be interesting to understand 
what happens when $M$ only has non-negative scalar curvature.  
Theorem \ref{poscal} still gives $[S]$ as an obstruction, but is very 
likely no longer sharp.  As noted in the sketch of the proof of theorem
\ref{constr}, the construction there does not give  metrics of bounded
geometry unless one has a positive lower bound on $\kappa$.  It seems 
likely that in the place of $l^{\infty}$-cohomology one needs forms 
which go to zero at $\infty$ in some way related to $\kappa$.  A closer 
examination of the behavior of $C$ in Theorem \ref{bndry} in terms of the 
lower bound on $\kappa$ might give the right decay condition.

   One intriguing aspect of Theorem \ref{poscal} is that the 
obstruction lives in a non-Hausdorff group.  This prevents the problem 
from fitting in to the $C^{*}$ algebra framework usually used for 
these types of problems.  To get around this, one can work with 
reduced uniformly finite homology, where the groups are the quotients 
of cycles by the closure of the boundaries.  The reduced invariant is 
shown to obstruct positive scalar curvature in \cite{Roe} and 
\cite{BW2}.  Corollary \ref{sums} shows that this loses important 
information.  The necessity of working with non-Hausdorff groups makes 
it unclear how to get versions of the theorems for families, and in 
particular, it is unclear when $(M \#_S N) \times {\mathbb R}^i$ carries 
a metric of positive scalar curvature.

  Theorem \ref{bndry} has applications to compact manifolds as well.  
Recall the theorem of \cite{Kaz} which says that for any compact 
manifold there is a metric whose non-positive scalar curvature is 
contained in an arbitrary ball.  This is a purely topological 
statement.  As a corollary of \ref{bndry}, one can see that this ball 
cannot be arbitrarily small.  

\begin{theorem} Let $M$ be a compact spin manifold with $\hat{A} \neq 
0$.  For any bounds on the curvature, there is an $r>0$ so that there 
is no metric on $M$ with those bounds on curvature and whose 
non-positive curvature set is contained in a ball of radius $r$.
\begin{proof}

  For $r$ small enough the ball is embedded.  Further, one has bounds 
on the second fundamental form of the sphere in terms of $r$ and the 
curvature of $M$.  Applying Theorem \ref{bndry} to the complement of 
the ball gives a bound, which goes to $0$ with $r$, on the integral 
of the $\hat{A}$ class over the complement of the ball.  As the 
curvature is bounded, the integral over the ball is bounded by a 
multiple of its volume.  These must add to a non-zero integer, which 
is a contradiction for $r$ sufficiently small.

\end{proof}
\end{theorem}

Also, many of the standard applications (see, for example, \cite{LM}) of the 
Bochner method can be carried over to manifolds with boundary as 
well: flat manifolds have their signature bounded by a multiple of the 
volume of the boundary, likewise for the total Betti number of 
manifolds with positive definite curvature tensor.  For submanifolds 
of Euclidean space or the sphere, these results follow easily just 
from Alexander duality.

\noindent
Kevin Whyte\\
Dept. of Mathematics\\
University of Chicago\\
5734 S. University Ave.\\
Chicago, Il 60637\\
E-mail: kwhyte@math.uchicago.edu

\end{document}